\documentstyle[12pt]{amsart}

\title{A Tverberg-type result on multicolored simplices}

\author{J\'anos Pach}

\thanks{Supported by NSF grant
CCR-94-24398,
PSC-CUNY Research Award 663472 and OTKA-4269.
This paper was written while the author was visiting
MSRI Berkeley, as part of the Convex 
Geometry Program. Research at MSRI is supported
in part 
by NSF grant DMS-9022140.}

\address{\hskip-\parindent J\'anos Pach\\
City College, CUNY and\\
Courant Institute, NYU\\
251 Mercer Street\\
New York, NY 10012}

\email{pach\@cims6.nyu.edu}

\date{}

\newtheorem{theorem}{Theorem}[section]

\def\bc{\begin{center}}\def\ce{\begin{center}}\def\ec{\end{center}}
\def\bd{\begin{eqnarray*}}\def\ed{\end{eqnarray*}}

\def\eps{\varepsilon}
\def\qed{\ifhmode\unskip\nobreak\fi\quad\ifmmode\Box\else$\Box$\fi}

\begin{document}

\begin{abstract}
Let $P_1, P_2,\ldots, P_{d+1}$ be pairwise disjoint 
$n$-element point sets in general position in 
$d$-space. It is shown that there exist
a point $O$ and suitable subsets $Q_i\subseteq P_i
\; (i=1, 2, \ldots, d+1)$ such that $|Q_i|\geq c_d|P_i|$,
and every $d$-dimensional 
simplex with exactly one vertex in each $Q_i$ contains
$O$ in its interior. Here $c_d$ is a positive constant
depending only on $d$.
\end{abstract}

\maketitle
\section{Introduction}

Let $P_1, P_2,\ldots, P_{d+1}$ be pairwise disjoint
$n$-element point sets in general position in Euclidean
$d$-space $\Re ^d$. If two points belong to the same $P_i$, then
we say that they are of the same {\it color}. A $d$-dimensional
simplex is called {\it multicolored}, if it has exactly
one vertex in each $P_i\; (i=1, 2, \ldots, d+1)$.
Answering a question of B\'ar\'any, F\"uredi, and 
Lov\'asz [BFL90], Vre\'cica and  {\v Z}ivaljevi\'c [ZV92],
proved the following Tverberg-type result. For every
$k$, there exists an integer $n(k,d)$ such that if 
$n\geq n(k,d)$, then one can always find $k$ multicolored 
vertex disjoint simplices with an interior point in common.
(For some special cases, see [BL92], [JS91], [VZ94].) This
theorem can be used to derive a nontrivial upper bound on the
number of different ways one can cut a finite point set
into two (roughly) equal halves by a hyperplane.

The aim of this note is to strengthen the above result
by showing that there exist ``large'' subsets of the sets
$P_i$ such that {\it all} multicolored simplices induced
by them have an interior point in common. 

\medskip
\medskip
\medskip
{\bf Theorem.} {\it There exists $c_d > 0$ with the property 
that for any disjoint $n$-element point sets
$P_1, P_2,\ldots, P_{d+1}\subset\Re ^d$ in general
position, one can find a point $O$ and suitable subsets
$Q_i\subseteq P_i,\; |Q_i|\geq c_d |P_i|\;
(i=1, 2, \ldots, d+1)$ such that every $d$-dimensional
simplex with exactly one vertex in each $Q_i$ contains
$O$ in its interior.}

\medskip 
The proof is based on the $k=d+1$ special case of the
Vre\'cica-\v Zivaljevi\'c theorem (see Theorem 2.1).
It uses three auxiliary results, each of them interesting 
on its own right. The first is Kalai's fractional Helly theorem
[K84], which sharpens and generalizes some earlier results
of Katchalski and Liu [KL79] (see Theorem 2.2).
The second is a variation of Szemer\'edi's regularity lemma
for hypergraphs [S78] (Theorem 2.3), and the third is a corollary 
of Radon's theorem [R21], discovered and applied by Goodman and 
Pollack [GPW96] (Theorem 2.4).

In the next section, we state the above mentioned results
and also include a short proof of Theorem 2.3, because in its 
present form it cannot be found in the literature. Our argument is
an adaptation of the approach of Koml\'os and S\'os [KS96].
For some similar results, see [C91],[FR92],[KS95]. The proof 
of the Theorem is given in Section 3. It shows
that the statement is true for a constant $c_d>0$ whose value
is triple-exponentially decreasing in $d$. 

\section{Auxiliary results}

\begin{theorem}\label{1}
{\rm [ZV92]} Let $A_1, A_2,\ldots, A_{d+1}$ be disjoint
$4d$-element sets in general position in $d$-space. Then one
can find $d+1$ vertex disjoint simplices with a common interior
point such that each of them has exactly one vertex in every
$A_i,\; 1\leq i\leq d+1$.
\end{theorem}

A family of sets is called {\em intersecting} if they have an
element in common.

\begin{theorem}\label{2}
{\rm [K84]} For any $\alpha > 0$, there exists
$\beta = \beta(\alpha,d) > 0$ satisfying the following
condition. Any family of $N$ convex sets in $d$-space,
which contains at least $\alpha {N\choose {d+1}}$ 
intersecting $(d+1)$-tuples, has an intersecting
subfamily with at least $\beta N$ members.
\end{theorem}

In fact, if $N$ is sufficiently large, then Theorem 2.2
is true for any $\beta < 1 - (1-\alpha)^{1/(d+1)}$. In
particular, it holds for $\beta = \alpha/(d+1)$.

Let $\cal H$ be a $(d+1)$-partite hypergraph whose vertex
set is the union of $d+1$ pairwise disjoint $n$-element
sets, $P_1, P_2,\ldots, P_{d+1}$, and whose edges are 
$(d+1)$-tuples containing precisely one element from each
$P_i$. For any subsets $S_i\subseteq P_i\; (1\leq i\leq d+1)$,
let $e(S_1, \ldots, S_{d+1})$ denote the number of edges of
$\cal H$ induced by $S_1\cup\ldots\cup S_{d+1}$. In this
notation, the total number of edges of $\cal H$
is equal to $e(P_1, \ldots, P_{d+1})$.

It is not hard to see that for any sets $S_i$ and for any
integers $t_i\leq |S_i|,\; 1\leq i\leq d+1$,

\begin{equation}\label{average}
\frac{e(S_1, \ldots, S_{d+1})}{|S_1|\ldots |S_{d+1}|}
=\sum\frac{e(T_1, \ldots, T_{d+1})}{|T_1|\ldots |T_{d+1}|}
/{{|S_1|}\choose {t_1}}\ldots {{|S_{d+1}|}\choose {t_{d+1}}},
\end{equation}

\noindent where the sum is taken over all $t_i$-element subsets
$T_i\subseteq S_i,\; 1\leq i\leq d+1$.

\begin{theorem}\label{3}
Let $\cal H$ be a $(d+1)$-partite hypergraph on the vertex
set $P_1\cup\ldots\cup P_{d+1},\; |P_i| = n\; (1\leq i
\leq d+1)$, and assume that $\cal H$ has at least
$\beta n^{d+1}$ edges for some $\beta > 0$. Let
$0<\eps <1/2$.

Then there exist subsets $S_i\subseteq P_i$ of equal
size $|S_i|=s\geq
\beta ^{1/\eps ^{2d}}n\; (1\leq i\leq d+1)$ such that

\noindent (i) $e(S_1, \ldots, S_{d+1})
\geq\beta s^{d+1},$

\noindent (ii) $e(Q_1, \ldots, Q_{d+1})>0$
for any $Q_i\subseteq S_i$ with
$|Q_i|\geq\eps s\; (1\leq i\leq d+1)$.
\end{theorem}

\noindent {\bf Proof:} Let $S_i\subseteq P_i\; (1\leq i
\leq d+1)$ be sets of equal size such that
$$\frac{e(S_1,\ldots, S_{d+1})}{|S_1|^{d+1-\eps ^{2d}}}$$
is {\em maximum}, and denote $|S_1|=\ldots =|S_{d+1}|$ 
by $s$.

For this choice of $S_i$, condition (i) in the theorem
is obviously satisfied, because
$$\frac{e(S_1,\ldots, S_{d+1})}{|S_1|^{d+1-\eps ^{2d}}}
\geq\frac{e(P_1,\ldots, P_{d+1})}{n^{d+1-\eps ^{2d}}}
=\frac{\beta}{n^{-\eps ^{2d}}}\geq
\frac{\beta}{s^{-\eps ^{2d}}}.$$
Taking into account the trivial relation
$$\frac{e(S_1,\ldots, S_{d+1})}{|S_1|^{d+1-\eps ^{2d}}}
\leq s^{\eps ^{2d}},$$
the above inequalities also yield that 
$s\geq\beta ^{1/\eps ^{2d}}n$.

It remains to verify (ii). To simplify the notation,
assume that $\eps s$ is an integer, and let $Q_i$ be
any $\eps s$-element subset of $S_i\; (1\leq i\leq d+1).$
Then

\begin{eqnarray*}
\lefteqn{e(Q_1,\ldots, Q_{d+1}) =  e(S_1,\ldots, S_{d+1})} \\
& & - e(S_1-Q_1, S_2, S_3, ..., S_{d+1})\\
& & - e(Q_1, S_2-Q_2, S_3, ..., S_{d+1})\\
& & - e(Q_1, Q_2, S_3-Q_3, ..., S_{d+1})\\
& & \ldots\\
& & - e(Q_1, Q_2, Q_3, ..., S_{d+1}-Q_{d+1}).
\end{eqnarray*}

In view of (1), it follows from the maximal choice of 
$S_i$ that

\begin{eqnarray*}
\lefteqn{e(S_1-Q_1, S_2, ..., S_{d+1})}\\
& & = (1-\eps)s^{d+1}\frac{e(S_1-Q_1, S_2, ..., S_{d+1})}
{|S_1-Q_1| |S_2| ... |S_{d+1}|}\\
& & = (1-\eps)s^{d+1}\sum_{{T_i\subseteq S_i, |T_i|=
(1-\eps)s\atop 2\leq i\leq d+1}}
\frac{e(S_1-Q_1, T_2, ..., T_{d+1})}{[(1-\eps)s]^{d+1}}
/{s\choose\eps s}^d\\
& & \leq(1-\eps)s^{d+1}\frac{e(S_1, S_2, ..., S_{d+1})}
{s^{d+1-\eps ^{2d}}}[(1-\eps)s]^{-\eps ^{2d}}\\
& & = e(S_1, ..., S_{d+1})(1-\eps)^{1-\eps ^{2d}}.
\end{eqnarray*}

\noindent Similarly, for any $i,\; 2\leq i\leq d+1$, we have

\begin{eqnarray*}
\lefteqn{e(Q_1, ..., Q_{i-1}, S_i-Q_i, S_{i+1}, ..., S_{d+1})}\\
& & \leq e(S_1, ..., S_{d+1})\eps ^{i-1-\eps ^{2d}}(1-\eps).
\end{eqnarray*}

Summing up these inequalities, we obtain

\begin{eqnarray*}
\lefteqn{e(Q_1,\ldots, Q_{d+1})\geq e(S_1,\ldots, S_{d+1})
(1-(1-\eps)^{1-\eps ^{2d}}-\sum_{i=2}^{d+1}
\eps ^{i-1-\eps ^{2d}}(1-\eps))}\\
& & \geq e(S_1,\ldots, S_{d+1})(1-(1-\eps)^{1-\eps ^{2d}}-
\eps ^{1-\eps ^{2d}}+\eps ^{d+1-\eps ^{2d}}) > 0,
\end{eqnarray*}

\noindent as required. \qed

A $(d+1)$-tuple of convex sets in $d$-space is called
{\em separated} if any $j$ of them can be strictly separated
from the remaining $d+1-j$ by a hyperplane, $1\leq j\leq d$.
An arbitrary family of at least $d+1$ convex sets in $d$-space is
{\em separated} if every $(d+1)$-tuple of it is separated.

\begin{theorem}\label{4}
{\rm [GPW96]} A family of convex sets in $d$-space is separated
if and only if no $d+1$ of its members can be intersected by
a hyperplane.
\end{theorem}

Let $n\geq d+1$. Two sequences of points in $d$-space,
$(p_1,\ldots, p_n)$ and $(q_1,\ldots, q_n)$,
are said to have the same {\em order type} if for any
integers $1\leq i_1 <\ldots <i_{d+1}\leq n$, the simplices
$p_{i_1}\ldots p_{i_{d+1}}$ and $q_{i_1}\ldots q_{i_{d+1}}$
have the same orientation [GP93]. It readily follows from the last
result that if $C_1,\ldots, C_n$ form a separated family of
convex sets, then the order type of $(p_1,\ldots, p_n)$ will
be the same for every choice of elements $p_i\in C_i,\;
1\leq i\leq n.$ 

\section{Proof of Theorem}

Let $P_1,\ldots, P_{d+1}$ be pairwise disjoint $n$-element point
sets in general position in $d$-space. If a simplex has precisely
one vertex in each $P_i$, we call it {\em multicolored}. The
number of multicolored simplices is $N=n^{d+1}$.

By Theorem 2.1, any collection of $4d$-element subsets $A_i\subseteq P_i,\;
1\leq i\leq d+1,$ induce $d+1$ vertex disjoint multicolored 
simplices with a common interior point. Thus, the total number
of intersecting $(d+1)$-tuples of multicolored simplices
is at least
$$\frac{{n\choose 4d}^{d+1}}{{{n-d-1}\choose {3d-1}}^{d+1}}>
\frac{1}{(5d)^{d^2}}{N\choose {d+1}}.$$
Hence, we can apply Theorem 2.2 with $\alpha = 1/(5d)^{d^2}$.
We obtain that there is a point $O$ contained in the interior
of at least $\beta N = \beta (1/(5d)^{d^2},d)n^{d+1}$ 
multicolored simplices.

Let $\cal H$ denote the $(d+1)$-partite hypergraph on the
vertex set $P_1\cup\ldots\cup P_{d+1}$, whose edge set 
consists of all multicolored $(d+1)$-tuples that induce a
simplex containing $O$ in its interior. 

Set $\eps = 1/2^{d2^d},$ and apply Theorem 2.3 to the
hypergraph $\cal H$ to find $S_i\subseteq P_i,\;
1\leq i\leq d+1,$ meeting the requirements. By throwing
out some points from each $S_i$, but retaining a positive
proportion of them, we can achieve that the convex hulls
of the sets $S_i$ are separated. Indeed, assume e.g. that there is
no hyperplane strictly separating $S_1\cup\ldots\cup S_j$
from $S_{j+1}\cup\ldots\cup S_{d+1}$. By the {\em ham-sandwich
theorem} [B33], one can find a hyperplane $h$ which 
simultaneously bisects $S_1,\ldots, S_d$ into as
equal parts as possible. Assume without loss of generality 
that at least half of the elements of $S_{d+1}$ are
``above'' $h$. Then throw away all elements of
$S_1\cup\ldots\cup S_j$ that are above $h$ and all elements
of $S_{j+1}\cup\ldots\cup S_{d+1}$ that are below $h$. 
We can repeat this procedure as long as we find a 
non-separated $(d+1)$-tuple. In each step, we reduce the 
size of every set by a factor of at most 2.

Notice that in the same manner we can also achieve that 
e.g. the $(d+1)$-tuple \(\{\{O\}, \mbox{ conv$(S_1),\ldots,$ 
conv$(S_d)$}\}\) becomes separated. In this case,  
$h$ will always pass through the point $O$,
therefore $O$ will never be deleted.

After at most $(d+2)2^d$ steps we end up with
$Q_i\subseteq S_i,\; |Q_i|>\eps s\; (1\leq i\leq d+1)$
such that \(\{\{O\}, \mbox{ conv$(S_1),\ldots,$conv$(S_{d+1})$}\}\)
is a separated family. It follows from the remark after
Theorem 2.4 that there are only two possibilities: either
every multicolored simplex induced by $Q_1\cup\ldots\cup Q_{d+1}$
contains $O$ in its interior, or none of them does. However,
this latter option is ruled out by part (ii) of Theorem 2.3.
This completes proof. \qed

\medskip

Instead of applying Theorem 2.2, we could have started the 
proof by referring to the following result of Alon,
B\'ar\'any, F\"uredi, and Kleitman [ABFK92], which is
also based on Theorem 2.1. For any
$\beta>0$ there is a $\beta^{'}_{d} > 0$ such that any family
of $\beta n^{d+1}$ simplices induced by $n$ points in
$d$-space has at least $\beta^{'}_{d}n^{d+1}$ members with
non-empty intersection.

Our proof easily yields the following. 

\begin{theorem}
For any $\beta>0$ there is a $\beta^{''}_{d}> 0$ with
the property that given any family of $\beta n^{d+1}$
simplices induced by an $n$-element set $P\subset\Re ^d$, 
one can find a point $O$ and suitable subsets
$Q_i\subseteq P_i,\; |Q_i|\geq\beta^{''}_{d}|P_i|\;
(i=1, 2, \ldots, d+1)$ such that every $d$-dimensional
simplex with exactly one vertex in each $Q_i$ contains
$O$.
\end{theorem} 

{\bf Acknowledgement.} I am grateful to Imre
B\'ar\'any, G\'eza T\'oth, and Pavel Valtr for their
valuable suggestions.

\end{document}